\documentclass[12pt]{amsart}
\usepackage{amsmath, amssymb, latexsym, amsthm, mathrsfs, color, mathtools, url}
\usepackage{extsizes, blindtext} 

\usepackage[top=2in, bottom=1.5in, left=1in, right=1in]{geometry}

\newcommand{\nc}{\newcommand}
\nc{\ea}{Markov} 
\nc{\lei}{\le^\oo}
\nc{\zar}{\mathfrak{zar}}
\nc{\seq}[1]{(#1)_{n\in\bbN}}
\nc{\card}[1]{\left|#1\right|}
\nc{\bbN}{\mathbb{N}}
\nc{\beq}{\begin{eqnarray*}}\nc{\eeq}{\end{eqnarray*}}
\nc{\mbq}{\mb{?}}
\nc{\mb}[1]{{\mbox{\textbf{#1}}}}
\nc{\nop}{$\times$}
\nc{\fbn}{\!\!\fbox{\!\nop\!}\!\!}
\nc{\yup}{\checkmark}
\nc{\forces}{\Vdash}
\nc{\name}[1]{\dot{#1}}
\nc{\tf}{\my{FINISHED THUS FAR}}
\nc{\FU}{Fr\'echet--Urysohn}
\nc{\gs}{$\gamma$~space}
\nc{\Ga}{\Gamma}\nc{\Om}{\Omega}
\nc{\smallbinom}[2]{\begin{psmallmatrix} #1\\ #2 \end{psmallmatrix}}
\nc{\bgamma}{\smallbinom{\Om}{\Ga}}
\newcommand{\two}{\{0,1\}}
\nc{\productive}[2]{\bigl(#1,\allowbreak #2\bigr)^\x}
\nc{\Sel}{\mathsf{S}}
\nc{\sset}[2]{{\{\,#1 : #2\,\}}}
\nc{\smb}[1]{{\!\!\mb{#1}\!\!}}
\nc{\cZ}{\mathcal{Z}}
\nc{\medset}[2]{{\biggl\{\,#1 : #2\,\biggr\}}}
\nc{\smallmedset}[2]{{\bigl\{\,#1 : #2\,\bigr\}}}
\nc{\set}[2]{{\left\{\,#1 : #2\,\right\}}}
\nc{\cube}{(\Cantor)^\bbN}
\nc{\Match}{\op{Match}}
\nc{\concat}[1]{\hat{\phantom{a}}\langle #1\rangle}
\nc{\poset}{\mathbb{P}}
\nc{\fn}[1]{{\op{Fn}(#1\times\w,2)}}
\nc{\linadd}{\op{linadd}}
\nc{\nonprod}{\non^\x}
\nc{\alephes}{{\aleph_0}}
\nc{\my}[1]{{\color{red} #1}}
\nc{\Cp}{\op{C}_\mathrm{p}}
\nc{\Bp}{\op{B}_p}
\nc{\Pa}[8]{\bibitem{#1} {#2}, \emph{#3}, {#4} \textbf{#5} ({#6}), {#7}--{#8}.}
\nc{\tPa}[5]{\bibitem{#1} {#2}, \emph{#3}, {#4}, to appear.}
\nc{\sPa}[4]{\bibitem{#1} {#2}, \emph{#3}, {#4}, submitted.}
\nc{\Bc}[9]{\bibitem{#1} {#2}, \emph{#3}, in: \textbf{#4} (#5), #6 #7, #8--#9.}
\nc{\fD}{\mathfrak{D}}
\nc{\fX}{\mathfrak{X}}
\nc{\Onbd}{\Op_{\mathrm{nbd}}} 
\nc{\Omnb}{\Om_{\mathrm{nbd}}} 
\nc{\od}{\mathfrak{od}}
\nc{\Setting}[7]{\xymatrix@R=4pt@C=7pt{#1\ar@{-}[r]&#2\ar@{-}[r]&#3\\&#4\ar@{-}[u]\\
#5\ar@{-}[uu]\ar@{-}[r] & #6\ar@{-}[u]\ar@{-}[r] & #7\ar@{-}[uu]}}
\nc{\mx}[1]{\begin{matrix}#1\end{matrix}}
\nc{\plim}{p\txt{-}\lim}
\nc{\Bgp}{{\Z^\bbN}}
\nc{\Cgp}{{{\Z_2}^\bbN}}
\nc{\Cite}[1]{\textbf{[#1]}}
\nc{\Next}[1]{{#1^+}}
\nc{\Fr}{\mathit{F\!r}}
\nc{\intvl}[2]{{[#1(#2),\allowbreak #1(#2\!+\!1))}}
\nc{\Bdd}{\mathbf{B}}
\nc{\Dfin}{\mathfrak{D}_\mathrm{fin}}
\nc{\grbl}{{\mbox{\textit{\tiny gp}}}}
\nc{\bbP}{\mathbb{P}}
\nc{\BOfat}{\B_{\Om_{\mathrm{fat}}}}
\nc{\Bgood}{\B_{\mathrm{good}}}
\nc{\compactN}{\cl{\mathbb{N}}}
\nc{\blocks}[2]{\op{cl}_{#2}(#1)}
\nc{\blocksplus}[2]{\op{cl}^+_{#2}(#1)}
\nc{\arx}[1]{\texttt{http://arxiv.org/math/#1}}
\nc{\bq}{\begin{quote}}
\nc{\eq}{\end{quote}}
\nc{\cl}[1]{\overline{#1}}
\nc{\CH}{the Continuum Hypothesis}
\nc{\MA}{Martin's Axiom}
\nc{\Bfat}{\B_\mathrm{fat}}
\nc{\inv}{^{-1}}
\nc{\Cantor}{{\two^\bbN}}
\nc{\bP}{\mathbf{P}}
\nc{\bof}{\op{\fb}}
\nc{\bofF}{\bof(\cF)}
\nc{\sr}[3]{\underset{\mbox{#3}}{\mbox{#1}}}
\nc{\gp}{\binom{\Om}{\Ga}}
\nc{\gpsmall}{\mbox{$\gp$}}
\nc{\gig}{\gimel}
\nc{\gns}{\sone(\Om,\gig)}
\nc{\nsr}[2]{#1}
\nc{\Srg}{{\mathbb{S}}}
\nc{\Srgs}{{\mathbb{S}^*}}
\nc{\NN}{{\bbN^{\bbN}}}
\nc{\ZN}{{\Z^{\bbN}}}
\nc{\NNup}{{\bbN^{\uparrow\bbN}}}
\nc{\PN}{{P(\bbN)}}
\nc{\roth}{{[\bbN]^{\mbox{\tiny $\infty$}}}} 
\nc{\Fin}{[\bbN]^{\text{$<\!\!\infty$}}} 
\nc{\ici}{[\bbN]^{ \infty, \infty}}
\nc{\Inc}{{\compactN^{\uparrow\bbN}}}
\nc{\powInc}[1]{{\big(\Inc\big)^{#1}}}
\nc{\powFin}[1]{{\big(\Fin\big)^{#1}}}
\nc{\powPN}[1]{{\big(\PN\big)^{#1}}}
\nc{\NcompactN}{{\compactN^\bbN}}
\nc{\Uarrow}{\smash{\big\uparrow}}
\nc{\LE}{\preccurlyeq}
\nc{\GE}{\succcurlyeq}
\nc{\op}{\operatorname}
\nc{\im}{\op{im}}
\nc{\Span}{\op{span}}
\nc{\maxfin}{\op{maxfin}}
\nc{\ran}{\op{range}}
\nc{\iso}{\cong}
\nc{\Madd}{{\M}^\star}
\nc{\cI}{\mathcal{I}}
\nc{\cJ}{\mathcal{J}}
\nc{\scrA}{\mathscr{A}}
\nc{\scrB}{\mathscr{B}}
\nc{\scrC}{\mathscr{C}}
\nc{\scrD}{\mathscr{D}}
\nc{\scrF}{\mathscr{F}}
\nc{\scrK}{\mathscr{K}}
\nc{\A}{\forall}
\nc{\B}{\mathrm{B}}
\nc{\cB}{\mathcal{B}}
\nc{\bB}{\mathbf{B}}
\nc{\BS}{\mathbf{B}(\mathcal{S})}
\nc{\BF}{\mathbf{B}(\mathcal{F})}
\nc{\BU}{\mathbf{B}(\mathcal{U})}
\nc{\cSp}{\mathcal{S}^+}
\nc{\cFp}{\mathcal{F}^+}
\nc{\cUp}{\mathcal{U}^+}

\nc{\BG}{\B_\Ga}
\nc{\BL}{\B_\Lambda}
\nc{\BT}{\B_\Tau}
\nc{\BTstar}{\B_{\Tau^*}}
\nc{\BO}{\B_\Om}
\nc{\DO}{\cD_\Om}
\nc{\KO}{\cK_\Om}
\nc{\CG}{C_\Ga}
\nc{\CL}{C_\Lambda}
\nc{\CT}{C_\Tau}
\nc{\CTstar}{C_{\Tau^*}}
\nc{\CO}{C_\Om}
\nc{\COgp}{C_{\Om^{\grbl}}}
\nc{\CLgp}{C_{\Lambda^{\grbl}}}
\nc{\BOgp}{\B_{\Om}^{\grbl}}
\nc{\BLgp}{\B_{\Lambda^{\grbl}}}
\nc{\sfC}{\mathsf{C}}
\nc{\sfD}{\mathsf{D}}
\nc{\bD}{\mathbf{D}}
\nc{\Tau}{\mathrm{T}}
\nc{\cA}{\mathcal{A}}
\nc{\cK}{\mathcal{K}}
\nc{\cD}{\mathcal{D}}
\nc{\cF}{\mathcal{F}}
\nc{\cS}{\mathcal{S}}
\nc{\cT}{\mathcal{T}}
\nc{\cG}{\mathcal{G}}
\nc{\cY}{\mathcal{Y}}
\nc{\J}{\mathcal{J}}
\nc{\cL}{\mathcal{L}}
\nc{\cM}{\mathcal{M}}
\nc{\cN}{\mathcal{N}}
\nc{\cE}{\mathcal{E}}
\nc{\cH}{\mathcal{H}}
\nc{\cO}{\mathcal{O}}
\nc{\Op}{\mathrm{O}}
\nc{\rmA}{\mathrm{A}}
\nc{\rmB}{\mathrm{B}}
\nc{\rmD}{\mathrm{D}}
\nc{\rmP}{\mathrm{P}}
\nc{\cC}{\mathcal{C}}
\nc{\cP}{\mathcal{P}}
\nc{\bbQ}{\mathbb{Q}}
\nc{\bbR}{\mathbb{R}}
\nc{\bbZ}{\mathbb{Z}}
\nc{\cU}{\mathcal{U}}
\nc{\Un}{\bigcup}
\nc{\cV}{\mathcal{V}}
\nc{\cW}{\mathcal{W}}
\nc{\Z}{{\mathbb Z}}
\nc{\Impl}{\Rightarrow}
\long\def\forget#1\forgotten{}
\nc{\ft}{\mathfrak{t}}
\nc{\fb}{\mathfrak{b}}
\nc{\fc}{\mathfrak{c}}
\nc{\fd}{\mathfrak{d}}
\nc{\fg}{\mathfrak{g}}
\nc{\oo}{\infty}
\nc{\fr}{\mathfrak{r}}
\nc{\fk}{\mathfrak{k}}
\nc{\bidi}{\mathfrak{bidi}}
\nc{\fu}{\mathfrak{u}}
\nc{\fh}{\mathfrak{h}}
\nc{\fp}{\mathfrak{p}}
\nc{\fj}{\mathfrak{j}}
\nc{\fs}{\mathfrak{s}}
\nc{\w}{\omega}
\nc{\x}{\times}
\nc{\Iff}{\Leftrightarrow}

\nc{\nin}{\notin}
\nc{\cat}{\hat{\ }}
\nc{\sub}{\subseteq}
\nc{\spst}{\supseteq}
\nc{\sm}{\setminus}
\nc{\as}{\subseteq^*}
\nc{\les}{\le^*}
\nc{\leinf}{\le^{\infty}}
\nc{\leS}{\le_{\mathcal{S}}}
\nc{\leF}{\le_{\mathcal{F}}}
\nc{\leU}{\le_{\mathcal{U}}}
\nc{\rest}{\restriction}

\nc{\la}{\langle}
\nc{\ra}{\rangle}
\nc{\E}{\exists}
\nc{\dom}{\op{dom}}
\nc{\cov}{\op{cov}}
\nc{\add}{\op{add}}
\nc{\cof}{\op{cof}}
\nc{\cf}{\op{cf}}
\nc{\non}{\op{non}}
\nc{\unif}{\op{non}}
\nc{\COV}{\op{COV}}
\nc{\ADD}{\op{ADD}}
\nc{\COF}{\op{COF}}
\nc{\NON}{\op{NON}}
\nc{\impl}{\to}
\nc{\Lp}{\mathcal{L_\p}}
\nc{\Wlog}{without loss of generality}
\newtheorem{thm}{Theorem}
\nc{\bthm}{\begin{thm}} \nc{\ethm}{\end{thm}}
\newtheorem{prop}[thm]{Proposition}
\nc{\bprp}{\begin{prop}} \nc{\eprp}{\end{prop}}
\newtheorem{fact}[thm]{Fact}
\nc{\bfct}{\begin{fact}} \nc{\efct}{\end{fact}}
\newtheorem{prob}[thm]{Problem}
\nc{\bprb}{\begin{prob}} \nc{\eprb}{\end{prob}}
\newtheorem{lem}[thm]{Lemma}
\nc{\blem}{\begin{lem}} \nc{\elem}{\end{lem}}
\newtheorem{claim}[thm]{Claim}
\nc{\bclm}{\begin{claim}} \nc{\eclm}{\end{claim}}
\newtheorem{cor}[thm]{Corollary}
\nc{\bcor}{\begin{cor}} \nc{\ecor}{\end{cor}}
\newtheorem{conj}[thm]{Conjecture}
\nc{\bcnj}{\begin{conj}} \nc{\ecnj}{\end{conj}}
\theoremstyle{definition}
\newtheorem{defn}[thm]{Definition}
\nc{\bdfn}{\begin{defn}} \nc{\edfn}{\end{defn}}
\newtheorem{obs}[thm]{Observation}
\nc{\bobs}{\begin{obs}} \nc{\eobs}{\end{obs}}
\theoremstyle{remark}
\newtheorem{rem}[thm]{Remark}
\nc{\brem}{\begin{rem}} \nc{\erem}{\end{rem}}
\newtheorem{cnv}[thm]{Convention}
\nc{\bcnv}{\begin{cnv}} \nc{\ecnv}{\end{cnv}}
\newtheorem{exam}[thm]{Example}
\nc{\bexm}{\begin{exam}} \nc{\eexm}{\end{exam}}
\nc{\bpf}{\begin{proof}} \nc{\epf}{\end{proof}}
\nc{\be}{\begin{enumerate}}
\nc{\ee}{\end{enumerate}}
\nc{\bi}{\begin{itemize}}
\nc{\bimy}{\my{\begin{itemize}}
\nc{\eimy}{\end{itemize}}}
\nc{\itm}{\item}
\nc{\ei}{\end{itemize}}
\nc{\ed}{\end{document}}
\nc{\Subsection}[1]{\goodbreak\subsection*{#1}}
\nc{\ffm}{\mathfrak{fm}}
\nc{\fmar}{\mathfrak{mar}}
\nc{\fkum}{\mathfrak{kum}}
\forget
\forgotten

\nc{\sone}{\mathsf{S}_1}
\nc{\sfin}{\mathsf{S}_\mathrm{fin}}
\nc{\ufin}{\mathsf{U}_\mathrm{fin}}
\nc{\Split}{\mathsf{Split}}

\nc{\gone}{\mathsf{G}_1}    \nc{\gfin}{\mathsf{G}_\mathrm{fin}}

\title{The Haar Measure Problem}

\author[A. Prze\'zdziecki]{Adam J. Prze\'zdziecki}
\address{Adam Prze\'zdziecki, Warsaw University of Life Sciences---SGGW, Warsaw, Poland}
\email{adamp@mimuw.edu.pl}
\urladdr{http://www.mimuw.edu.pl/~adamp}

\author[P. Szewczak]{Piotr Szewczak}
\address{Piotr Szewczak, Faculty of Mathematics and Natural Science College of Sciences, Cardinal Stefan Wyszy\'nski University in Warsaw, Warsaw, Poland, and Department of Mathematics, Bar-Ilan University, Ramat Gan, Israel,}
\email{p.szewczak@wp.pl}
\urladdr{www.piotrszewczak.pl}

\author[B. Tsaban]{Boaz Tsaban}
\address{Boaz Tsaban, Department of Mathematics, Bar-Ilan University, Ramat Gan, Israel}
\email{tsaban@math.biu.ac.il}
\urladdr{http://math.biu.ac.il/~tsaban}

\subjclass[2010]{
	28C10,  	
	28A05,  	
	22C05,  	
	03E17
}
\keywords{Haar measurable, Baire property, profinite group, compact group, closed measure zero}

\begin{document}

\begin{abstract}
  An old problem asks whether every compact group has a Haar-nonmeasurable subgroup. A series of earlier results reduce the problem to infinite metrizable profinite groups. 
  We provide a positive answer, assuming a weak, potentially provable,
   consequence of \CH{}. We also establish the dual, Baire category analogue of this result.
\end{abstract}

\maketitle

\section{Introduction}
\label{sec:intro}

Every infinite compact group has a unique translation-invariant probability measure,
its \emph{Haar measure}. 
Vitali sets (complete sets of coset representatives) with respect to a
countably infinite subgroup show that such groups 
have nonmeasurable subsets. We consider the following old problem.

\theoremstyle{theorem}
\newtheorem*{hmp}{Haar Measure Problem}
\begin{hmp}
Does every infinite compact group have a nonmeasurable subgroup?
\end{hmp}

The Haar Measure Problem dates back at least to 1963, 
when Hewitt and Ross gave a positive answer for abelian groups~\cite[Section~16.13(d)]{HR}.
It was explicitly formulated in a paper of Saeki and Stromberg~\cite{SS85}.
The problem remains open despite substantial efforts~\cite[and references therein]{HHM, BM}.

 Hern\'andez, Hofmann, and Morris proved that if all subgroups of an infinite compact
 group are measurable, then the group must be profinite and metrizable~\cite[Theorem 2.3 and Corollary~3.3]{HHM}.
Building on that, Brian and Mislove proved that a positive answer to the Haar Measure
Problem is \emph{consistent} (relative to the usual axioms of set theory)~\cite[Theorem~2.5]{BM}.
We repeat their argument, for its elegant simplicity, and since this result will take
care of the easier case of our main theorem: 
Let $G$ be an infinite, metrizable profinite group. As a measure space,
the group $G$ is isomorphic to the Cantor space with the Lebesgue measure.
Let $\fc$ denote the cardinality of the continuum.
Consistently, 
there is in the Cantor space, and thus in $G$,
a nonnull set $A$ of cardinality smaller than $\fc$. 
The subgroup of $G$ generated by $A$ is nonnull, and its
cardinality is smaller than $\fc$.
Since sets of positive measure have cardinality
 $\fc$, the group generated by $A$ is nonmeasurable.

Brian and Mislove's observation can be viewed as a solution of the Haar Measure
Problem under the hypothesis that there is a nonnull set of cardinality smaller than
$\fc$. This hypothesis violates \CH{}. We will show that \CH{} also implies a positive solution. Moreover, for our proof we only assume a weak consequence of \CH{}, 
which is provable for some groups, and may turn out provable for all groups. 
A proof of our hypothesis, if found, would settle the Haar Measure
Problem.

\section{The main theorem}

Throughout this section, we fix an arbitrary infinite metrizable profinite group $G$, 
and let $\mu$ be its Haar probability measure. For each natural number $n$, the
Haar probability measure on the group $G^n$ is the product measure, which is also denoted
$\mu$.

Let $H$ be a subgroup of $G$, and  $X=\{x_1,x_2,\dotsc\}$ be a countable set of variables.
The set of all words in the alphabet $H\cup\{x_1^{\pm 1},x_2^{\pm 1},\dotsc\}$ is denoted $H[X]$. 
Each word $w\in H[X]$ depends on finitely many parameters from the subgroup $H$, and finitely many variables;
let $\card{w}$ denote the number of variables in $w$. 
We view the word $w$ as a continuous function from $G^{\card{w}}$ to $G$
defined by substituting the group elements for the variables.

\bdfn
Let $e$ be the identity element of the group $G$.
A \emph{\ea{} set} is a set of the form
$w\inv(e)$ for $w\in G[X]\sm G$.
\edfn

The \ea{} sets were studied by Markov, as the sets that
are closed in all group topologies on $G$. 

\blem
For each 
element $b\in G$, and each word $w\in G[X]\sm G$, the set $w\inv(b)$ is \ea{}.
\elem
\bpf
Consider the word $wb\inv$.
\epf

For a natural number $n\ge 2$, a set $A\sub G^n$, and an element $g\in G$, we define 
\[
A_g:=\set{h\in G^{n-1}}{(h,g)\in A},
\]
the fiber of the set $A$ over the point $g$ in the group $G^{n-1}$.

Since \ea{} sets are closed (and thus measurable), so are their fibers.

\bdfn\label{dfn:FM}
A \emph{\ea{} null set} is a \ea{} subset of some
finite power of the group $G$, that is also null with
respect to the Haar measure $\mu$.
A set $N$ is \emph{Fubini--Markov} if either of the following two cases holds:
\be
\item The set $N$ is a \ea{} null subset of $G$.
\item There are a natural number $n\ge 2$ and a \ea{} null set 
$A\sub G^n$ such that 
$N=\set{g\in G}{\mu(A_g)>0}$.
\ee
\edfn

While \ea{} sets may be subsets of an arbitrary power of $G$,
Fubini--Markov sets are always subsets of $G$.
By the Fubini Theorem, we have the following observation.

\blem
\label{lem:fub}
Every Fubini--Markov set is null.\qed
\elem

We define a cardinal invariant of the group $G$.

\bdfn\label{dfn:FMnumber}
The \emph{Fubini--Markov number} of $G$, denoted $\ffm(G)$, is the minimal number of Fubini--Markov sets in $G$ whose union has full measure.
\edfn

Since a countable union of null sets is null, the Fubini--Markov number of a group is
necessarily uncountable.

\bexm
\label{exm:abel}
For the Cantor group, we have $\ffm(\bbZ_2^\bbN)=\fc$.
Indeed, let $G$  be an \emph{abelian} infinite metrizable profinite group, and $N$
be a Fubini--Markov set. We consider the two cases in the definition.

Assume that $N=w\inv(0)\sub G$ for some one-variable word $w\in G[X]\sm G$.
Since the group $G$ is abelian, we have $w(x)=x+a$ for some $a\in G$.
Then $w(x)=0$ if and only if $x=-a$, and thus $N$ is a singleton.

Next, for a natural number $n\ge 2$, let $w\inv(0)\sub G^n$ be a \ea{} null 
set, where $w\in G[X]\sm G$ and $\card{w}=n$. Since the group $G$ is abelian,
we have
$w:=x_1+\dotsb+x_n+a$, for some $a\in G$.
Then, for each $g\in G$, we have 
\[
(w\inv(0))_g=\set{(h_1,\dotsc,h_{n-1})}{h_1+\dotsb+h_{n-1}+g+a=0},
\]
and thus $(w\inv(0))_g$ 
is a Lipschitz image of the null set $G^{n-2}\x\{0\}$.
It follows that $\mu((w\inv(0))_g)=0$ for all $g\in G$, and 
we have, in the definition, $N=\emptyset$.

A union of fewer than $\fc$ sets that are at most singletons cannot cover a full measure set.
\eexm

We arrive at our main theorem.
Let $\cN$ be the ideal of null sets in the Cantor space, and $\non(\cN)$ be the minimal cardinality of a nonnull subset of the Cantor space.
We settle the Haar Measure Problem for groups $G$ with $\non(\cN)\le\ffm(G)$.
By Lemma~\ref{lem:fub}, \CH{} implies $\non(\cN)=\ffm(G)$. 
Example~\ref{exm:abel} shows that for some groups our 
hypothesis is provable. 
We do not know whether it is provable for all infinite metrizable profinite
groups $G$. The numbers $\ffm(G)$ are provably larger than some
classical cardinal invariants of the continuum; we will return to this in Section~\ref{sec:hyp}.

\bthm
\label{thm:main}
Let $G$ be an infinite metrizable profinite group with $\non(\cN)\le\ffm(G)$. 
Then $G$ has a Haar-nonmeasurable subgroup.
\ethm
\bpf
If $\non(\cN)<\fc$, then the Brian--Mislove argument applies, namely, every 
nonnull set of cardinality $\non(\cN)$ generates a nonmeasurable subgroup of $G$;
see Section~\ref{sec:intro} for the details. 
We may thus assume that $\non(\cN)=\fc$. By our hypothesis, we have $\ffm(G)=\fc$.

Let $\sset{N_\alpha}{\alpha<\fc}$ be the family of $\text{G}_\delta$ null sets. Every null set
is contained in some $N_\alpha$.
We define a transfinite, increasing chain of subgroups $H_\alpha$ of $G$ for $\alpha<\fc$.
Let $H_0$ be a countable dense subgroup of $G$.
Let $w\in H_0[X]\sm H_0$. For distinct elements $b_1,b_2\in G$, the sets
$w\inv(b_1)$ and $w\inv(b_2)$ are disjoint. Since \ea{} sets are closed
(and thus measurable), the set 
\[
P_w:=\set{b\in G}{\mu(w\inv(b))>0}
\] 
is countable.
Since the set $H_0[X]\sm H_0$ is countable, there is an element
\[
b\in G\sm \Bigl(H_0\cup \Un_{w\in H_0[X]\sm H_0}P_w\Bigr).
\]
This element $b$ will remain outside our subgroups throughout the construction.

We proceed by induction.
For a limit ordinal $\alpha$, we set  $H_\alpha:=\Un_{\beta<\alpha}H_\beta$.
For a successor ordinal $\alpha=\beta+1<\fc$,
we assume, inductively, that $\card{H_\beta}<\fc$, $b\nin H_\beta$, and the sets $w\inv(b)$ are null for all words  $w\in H_\beta[X]\sm H_\beta$.

Since $\card{H_\beta[X]}<\fc$, the set
\[
S:=\Un_{\substack{w\in H_\beta[X]\\
		\card{w}=1}}w\inv(b)\cup \Un_{\substack{w\in H_\beta[X]\\ 
		\card{w}\geq 2}}\smallmedset{g\in G}{\mu\bigl( (w\inv(b))_g\bigr)>0}.
\]
is a union of fewer than $\ffm(G)$
Fubini--Markov sets, and thus does not have full measure.
Pick an element $g_\alpha\in G\sm(S\cup N_\alpha)$.
Let $H_\alpha:=\langle H_\beta,g_\alpha\rangle$.
We verify that the inductive hypotheses are preserved.

Fix $c\in H_\alpha$. There is a word $w\in H_\beta[X]$ with $\card{w}=1$ such that $w(g_\alpha)=c$, and thus $g_\alpha\in w\inv(c)$.
Since $g_\alpha\notin S\spst w\inv(b)$, we have $c\neq b$. This shows that $b\notin H_\alpha$.
Next, consider an arbitrary word $v=v(x_1,\dotsc,x_n)\in H_\alpha[X]\sm H_\alpha$.
There is a word
$w=w(x_1,\dotsc,x_n,x_{n+1})\in H_\beta[X]\sm H_\beta$ such that
\[
v(x_1,\dotsc,x_n)=w(x_1,\dotsc,x_n,g_\alpha).
\]
Since $g_\alpha\nin S$ and $\card{w}\geq 2$, the set $v\inv(b)=(w\inv(b))_{g_\alpha}$ is null.

Having defined all subgroups $H_\alpha$ for $\alpha<\fc$, let $H:=\Un_{\alpha<\fc}H_\alpha$.
Then $H$ is a proper (since $b\notin H$), dense (since $H_0\sub H$), nonnull
(since $H\nsubseteq N_\alpha$ for all $\alpha$) subgroup of $G$.
Assume that $H$ is measurable. Then it has positive measure, and
by the Steinhaus Theorem, it contains an open set.
Since it is dense, we have $H=G$, a contradiction.
\epf

Our main theorem also has a dual, Baire category version.
Let $\cM$ be the ideal of meager (Baire first category) subsets of the Cantor space.
We define \emph{Kuratowski--Ulam--Markov sets} by 
changing \emph{null} to \emph{meager}
in Definition~\ref{dfn:FM}. In this case, item~(2) 
of the definition becomes 
\[
N=\sset{g\in G}{A_g\text{ is nonmeager}}.
\]
By the Kuratowski--Ulam Theorem, Kuratowski--Ulam--Markov sets are meager.
Similarly, we dualize Definition~\ref{dfn:FMnumber} to define
the Kuratowski--Ulam--Markov number $\fkum(G)$. 
Let $\non(\cM)$ be the minimal cardinality of a nonmeager subset of the Cantor space.

\bthm\label{thm:mainkat}
Let $G$ be an infinite metrizable profinite group with $\non(\cM)\le\fkum(G)$. 
Then $G$ has a subgroup that does not have the property of Baire.
\ethm
\bpf
If $\non(\cM)<\fc$, then any nonmeager set of cardinality $\non(\cM)$ generates
a nonmeager subgroup of $G$ that does not have the Baire property (nonmeager
sets with the Baire property have cardinality $\fc$)~\cite[Theorem~2.5]{BM}.

Thus, assume that $\non(\cM)=\fc$. 
We proceed as in the proof of Theorem~\ref{thm:main}, replacing \emph{null} by \emph{meager}
and sets of positive measure by \emph{nonmeager sets} (the relevant sets are closed).
For the choice of the element $b$,
we observe that closed nonmeager sets are, in particular, not nowhere dense, and thus
have nonempty interior. Our group $G$ is homeomorphic to the Cantor space, 
and thus there are in $G$ at most countably many disjoint open sets.

We thus obtain a proper dense nonmeager subgroup of $G$.
To conclude the proof, we use the Pettis Theorem~\cite[Theorem~9.9]{Kech}, the 
category-theoretic dual of the
Steinhaus Theorem: If a set $H\sub G$ is nonmeager and has the Baire property, then 
the quotient $H\inv H$  
has nonempty interior.
\epf

\section{Bounds on the Fubini--Markov number}
\label{sec:hyp}

Here too, all groups are assumed to be infinite metrizable profinite.
Theorem~\ref{thm:main} applies to groups $G$ with $\non(\cN)\le\ffm(G)$.
We saw that abelian groups have $\ffm(G)=\fc$, but the following conjecture
remains open.

\bcnj
For each infinite metrizable profinite group $G$, we have $\non(\cN)\le\ffm(G)$.
\ecnj

A proof of this conjecture would settle the Haar Measure Problem, but it may
turn out unprovable (and thus undecidable). In this case, well-studied lower bounds on 
$\ffm(G)$ are useful.

The \emph{covering number} of an ideal $\cI$ of subsets of the Cantor space, 
denoted $\cov(\cI)$, is the minimal number of elements of $\cI$ needed to cover
the Cantor space.
Since Fubini--Markov sets are null, we have
$\cov(\cN)\le\ffm(G)$ for all groups $G$: A set of full measure needs just one
additional null set to cover the entire space. The following result provides a
tighter estimate, in the sense that it is provably larger, and consistently strictly larger.

Let $\cE$ be the $\sigma$-ideal generated by the closed null sets in the Cantor space.
The following proof establishes, in particular, that the family of Fubini--Markov sets
is contained in $\cE$.

\bprp
\label{prp:relations}
For each infinite metrizable profinite group $G$, we have $\cov(\cE)\leq \ffm(G)$.
\eprp
\bpf
Brian proved that $\cov(\cE)$ is equal to the minimal number of closed null subsets of the Cantor space that cover a set of positive measure~\cite{Bri}. 
Thus, it suffices to prove that
every Fubini--Markov subset $N$ of $G$ is a countable union of closed null subsets of $G$.

Let $N$ be a Fubini--Markov subset of $G$. 
If $N$ is \ea{} null, then it is closed and null.
It remains to consider the case that
\[
N=\smallmedset{g\in G}{\mu\bigl( (w\inv(e))_g\bigr)>0},
\]
where $w\in G[X]$ has $\card{w}\geq 2$. 

For each natural number $k$,
the subset
\[
N_k:=\sset{g\in G}{\mu\bigl( (w\inv(e))_g  \bigr)\geq 1/k}
\]
of $N$ is null (Lemma~\ref{lem:fub}), and $N=\Un_kN_k$.
Each set $N_k$ is closed:
Let $g\in G\sm N_k$. There is an open set $V$ in $G^{\card{w}-1}$ such that $(w\inv(e))_g\sub V$ and $\mu(V)<1/k$. Let $P$ be the projection of the compact set $w\inv(e)\sm (V\x G)$ on the last coordinate.
The set $G\sm P$
is an open neighborhood of $g$ in $G$. 
For each element $h\in G\sm P$, we have $(w\inv(e))_{h}\sub V$.
Thus, $\mu((w\inv(e))_{h})\leq\mu(V)< 1/k$, and $(G\sm P)\cap N_k=\emptyset$.
\epf

\bcor
Assume that $\non(\cN)\le\cov(\cE)$.
Then every infinite compact group has a nonmeasurable subgroup.
\ecor
\bpf
Theorem~\ref{thm:main} and Proposition~\ref{prp:relations}.
\epf

Since $\cE\subseteq \cM\cap\cN$, we have 
\[
\max\{\cov(\cM),\cov(\cN)\}\le\cov(\cE).
\]
It follows that if $\non(\cN)\le\max\{\cov(\cM),\cov(\cN)\}$, then
every compact group has a nonmeasurable subgroup.
The hypothesis $\non(\cN)\le\cov(\cE)$ is not provable; this follows from 
known upper bounds on $\cov(\cE)$~\cite{Ash}.

We conclude this section with a simple sufficient condition for our main theorem.
This condition is stronger than the hypothesis $\non(\cN)\le\ffm(G)$, but it may still be provable. 

\bdfn
\label{dfn:mar}
Let $G$ be an infinite metrizable profinite group.
For a natural number $n$, let $\kappa_n$ be the minimal number of \ea{} null subsets of the group $G^n$ whose union is not null. 
The \emph{Markov number} of $G$ is the cardinal number
$\fmar(G):=\min_n\kappa_n$.
\edfn

\bprb
In Definition~\ref{dfn:mar}, is the sequence $\kappa_1,\kappa_2,\dots$ constant?
In particular, is it provable that $\fmar(G)$ is equal to the
minimal number of \ea{} null subsets of the group $G$ whose union is not null?
\eprb

\blem\label{lem:relations2}
Let $G$ be an infinite metrizable profinite group.
Then:
\be
\item $\cov(\cM) \leq \fmar(G) \leq \non(\cN)$,
\item $\fmar(G)\leq \ffm(G)$.
\ee
\elem
\bpf
(1)
\ea{} sets are closed and null.
It follows that the minimal number of closed null sets in the Cantor space
whose union is not null is at most $\fmar(G)$.
The former number is equal to $\cov(\cM)$~\cite[Theorem~2.6.14]{BarJu}.
Since every singleton is a \ea{} set (consider the word $w(x)=x$), 
we have $\fmar(G)\leq\non(\cN)$.

(2)
Let $\cF$ be a family of Fubini--Markov subsets of the group $G$
with $\card{\cF}<\fmar(G)$.
For each element of $\cF$, fix a \ea{} null set witnessing 
its being Fubini--Markov, and let $\cA$ be the family of these \ea{} sets.
For a natural number $n$,
let $A_n:=\Un\sset{A\in\cA}{A\sub G^n}$.
Since $\card{\cA}<\fmar(G)$, the set  $A_n$ is null.
By the Fubini Theorem, the set
\[
S:=A_1\cup \Un_{n=2}^\infty
\smallmedset{g\in G}{ \bigl( A_n\bigr)_g \text{ is not null}}
\]
is null.
Then $S$ is null, and $\Un\cF\sub S$.
\epf

Thus,
the following conjecture implies a positive solution to the Haar Measure Problem.

\bcnj 
\label{cnj:kappa}
For each infinite metrizable profinite group $G$, we have $\non(\cN)=\fmar(G)$.
\ecnj

Conjecture \ref{cnj:kappa} holds when restricted to abelian groups, since 
Lipschitz images of sets of the form $A\x \{0\}$, with $\card{A}<\non(\cN)$, are
null (see Example~\ref{exm:abel})

\subsection*{Acknowledgments.}
We are indebted to Will Brian for answering a question of ours~\cite{Bri}.
His answer simplified Proposition~\ref{prp:relations}. We thank the editor and
the referee for the evaluation of this paper.

\ed